\documentclass[]{interact}

\usepackage{epstopdf}
\usepackage[caption=false]{subfig}

\usepackage[longnamesfirst,sort]{natbib}
\bibpunct[, ]{(}{)}{;}{a}{,}{,}
\renewcommand\bibfont{\fontsize{10}{12}\selectfont}

\theoremstyle{plain}
\newtheorem{theorem}{Theorem}[section]

\theoremstyle{definition}
\newtheorem{definition}[theorem]{Definition}

\theoremstyle{remark}

\begin{document}

\articletype{ARTICLE TEMPLATE}

\title{New model for the decision-making and the game style in football}

\author{
\name{Brahim Boudine\textsuperscript{a}\thanks{CONTACT B. Boudine. Email: brahimboudine.bb@gmail.com} }
\affil{\textsuperscript{a} Faculty of sciences Dhar El Mahraz, Sidi Mohamed Ben Abdellah University, Fez, Morocco}
}

\maketitle

\begin{abstract}
We propose a new mathematical model for the decision-making of players in football (soccer) and the efficiency of the game style. Our approach is based on $4$-networks, which is a mathematical concept that we introduce. The decision of players is expressed by a mathematical function depending on the game style chosen by the coach. Moreover, we measure the efficiency of the game style by a sequence of $4$-networks. 
\end{abstract}


\section{Introduction}
Football is not just a game; it means so much more to our society. It brings people together, crossing borders, cultures, and backgrounds. The sport teaches valuable lessons like teamwork, discipline, and perseverance, empowering individuals to face life's challenges head-on.  For example, the king of Morocco Mohammed VI, in his speech on July 29, 2023, takes the example of the Moroccan national team in the world cup Qatar 2022 to inspire his people to be more confident and serious. In fact, the impact of football is significant in various fields. In 2007, Aza Conejo et al. \cite{Aza} studied its economic impact on the regional economy. In 2014, sociological impacts have been studied by Fillis and Mackay \cite{Fillis}, as well as, Spaaij \cite{Spaaij} investigated some political impacts of football. For these reasons, it is natural that the study of the ways of playing football becomes more interesting for scientists.

Actually, mathematics plays a crucial role in providing rigorous analysis across various fields. It offers a structured approach to problem-solving, distilling complex issues into precise equations for deeper understanding and accurate predictions. Therefore, it was important to provide mathematical models for analyzing football strategies. In 1997, with a statistical approach, Pollard and Reep \cite{Pollard} introduced a variable called yield to quantify and compare the effectiveness of strategies during possession. Since then, many authors was interested by the mathematics of football, we refer especially to Hughes \cite{Hughes} who provided a mathematical perspective, and  Szczepa?ski \cite{Szcz} who measures the effectiveness of strategies and quantifies players' performance. In the last decade, football continues to attract mathematicians, and mathematical analysis becomes more interesting for players and coaches.  In 2015, Gonzalez-Rodenas et al. \cite{Gonzalez3} investigated the association between tactics and scoring opportunities. In 2017, Dufour et al. \cite{Dufour} analyzed the world cup Brazil 2014 to determine the parameters which made the difference. Recently, many authors are interested in evaluating the tactical and technical performances of players, we refer to Wang et al. \cite{Wang}, Liu et al. \cite{Liu}, and Gonzalez-Rodenas et al. \cite{Gonzalez2}.

In this paper, we use networks to study how a player should make his decision to pass the ball or to shoot it. This is important for enhancing the decision-making of players, and providing a good evaluation plan to coaches. A player's decision is crucial to apply tactics, so it depends on tactical strategies given by the coach, and on the technical performance of players. Then, we introduce the notion of $n$-networks, which generalizes the known networks, where $n$ is a positive integer. For every couple of players $(t_i,t_j)$, we associate an $n$-edge, which allows us to show different parameters that we need to measure the function of the decision for the player $t_i$. This allows to transform tactics into mathematical functions. Afterwards, a sequence of $4$-networks shows the efficiency of the sequence of taken decisions; namely, the efficiency of the coach's tactic.   

\section{Background}
Graphs and networks have found extensive applications in various fields. Notably, in team sports like football, where players' decision-making plays a pivotal role in executing the coach's strategies and achieving efficient teamwork, these tools are highly beneficial. By employing graphs and networks, football teams can enhance their decision-making processes, optimize tactical approaches, and foster efficient collaboration among players. In 2015, Trequattrini et al. \cite{Trequattrini} provided the first application of network analysis of football team performance. Then, Clemente et al. \cite{Clemente} proposed some network methods to analyse the team's performance. By utilizing their approach, it becomes possible to discern the interactions between players, as well as the nature and intensity of their connections with one another. Afterwards, in another work, Clemente et al. \cite{Clemente2} studied the role of midfielders in the building attack, and they analyzed their positions using networks.

In modern football, we can notice that passing becomes a key skill for players. Therefore, particular importance has been given recently to analysing and evaluating passing ability. In 2016, with a statistical model, Szczepański and McHale \cite{Szcz2} evaluated passing ability and discerned some factors which affect the quality of the pass like the skill of the player and the defensive pressure. In 2018, Buldi et al. \cite{Buldu} used networks theory to analyze football passing networks, which are visual representation that illustrates a team's playing patterns based on their passing actions.  

Since a good execution of coaches' strategies involves good decision-making from players, it is important to determine the best decision to make by players. The above statistical works have already provided some models to measure different parameters like scoring ability, passing accuracy and passing quality. But, we need to measure how the decision taken is good in terms of tactics and game styles. This is important to evaluate the adaptation of the players with the coach's thoughts and the efficiency of the game style.

\section{Construction of our model}
In this section, we suppose that our team possess the ball. Let $T = \{t_1, t_2,...,t_{11}\}$ be the team; namely, each $t_i$ presents a player. Assume that the ball is with the player $t_i$, for a certain $i \in \{1,2,...,11\}$. Then, we consider the graph $G = (T,E_i)$, where $E_i = \{ \{t_i, t_j\} \mid j \neq i \emph{ and } 1\leq j \leq 11\}$ is the set of edges. So we put the questions:
\begin{itemize}
\item What should the player $t_i$ do with the ball?
\item How could he make the best decision? 
\end{itemize}  
In general, the player $t_i$ has two choices: either he shoots the ball or he passes it. Actually, this depends on his feelings: if he feels that he has a good chance to score, then he shoots the ball. Else, it is preferable to pass it. So, the first parameter that we should take is the probability of scoring $s_i \in [0,1]$. In 2014, McHale and Szczepański \cite{McHale} proposed a model to identify goal scoring ability of players. However, the ability of players is not sufficient to measure the probability of scoring. Recently, more works have been interesting in enhancing measuring the probability of scoring. In 2020, Gonzalez-Rodenas et al. \cite{Gonzalez} determined more indicators related to goal scoring, especially, technical, tactical and spacial indicators. Later, in 2021, Anzer and Bauer \cite{Anzer} measured the goal-scoring probability using a model based on the position of the player and the quality of his shot.  
\begin{figure}[h!]
\centering
{%
\resizebox*{15cm}{!}{\includegraphics{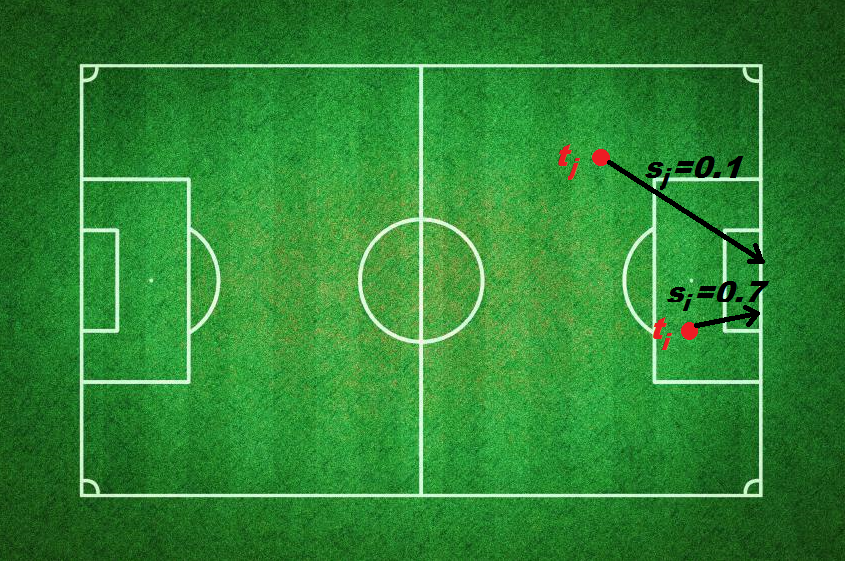}}}\hspace{3pt}
\caption{Example of a situation where two players with different probability of goal-scoring} \label{figure1}
\end{figure}

In figure \ref{figure1} the player $t_i$ has a good chance to score with a high probability $s_i=0.8$, so he can shoot the ball, but the player $t_j$ has a low probability to score $s_j=0.1$, so it is preferable to pass the ball. As well, if the player should pass the ball, then we put another question:
\begin{itemize}
\item How could he choose the teammate to whom he passes the ball?
\end{itemize}
Now, we consider three parameters to make this decision: 
\begin{enumerate}
\item The time $\tau_i \in \mathbb{R}_+$ in which the player $t_i$ should apply his decision. Notice that the time $\tau_i$ is low when the player $t_i$ is in high pressure, so he should make a fast decision. In this case, some qualities may play a significant role, such as control techniques and dribbling \cite{Yi}. Also, reading teammates' movements before receiving the ball (scanning skills), studied by Aksum et al. \cite{Aksum} in 2021, is important to optimize the necessary time for decision-making.   
\item The probability $p_{i,j}(\tau_i) \in [0,1]$ that the player $t_j$ will receive the ball from the player $t_i$ in comfortable conditions. This probability depends on the time $\tau_i$ in which the player $t_i$ should make his decision. Indeed, if the player $t_j$ is covered or strongly pressed by an adversary's defender, then the probability is low. The distance between the players $t_i$ and $t_j$, and the technical performance of the player $t_i$ will affect also the probability $p_{i,j}$. To simplify the notation, we will denote it only $p_{i,j}$. There are already some attempts to measure this metric, for example, we refer to \cite{Buldu} and \cite{Szcz2}.
\item The risk $r_j \in \{0,1,2,...,10\}$ which measures the threatening of the adversary's goal if the ball reaches the player $t_j$. This depends on the technical qualities of the player $t_j$, his probability to score, and his comfortability to make a better decision than the player $t_i$ who may be under high pressure.  
\end{enumerate}
If a player $t_j$ is offside or outside, then let us take $p_{i,j}=r_j=0$. \\

So, we have $4$ considerable parameters to make a decision. Then, we need a new mathematical concept that we call a $4$-network.
\begin{definition} Let $n$ be a positive integer. An $n$-network is a graph $G$ where every edge $e_i$ is equipped with an $n$-vector $v_i = \begin{pmatrix}
x_{i,1} \\
x_{i,2} \\
... \\
x_{i,n}
\end{pmatrix} \in \mathbb{R}^n$, the edge $e_i$ is called an $n$-edge and denoted $(e_i,v_i)$. As well, if $E$ is the set of edges of $G$, then we denote a $n$-network $N$ as follows $N=\{(e_i,v_i) \mid e_i \in E \emph{ and } v_i \in \mathbb{R}^n\}.$ 
\end{definition} 
For every $i,j \in \{1,2,...,11\}$, using the above notations, let $v_{i,j} = \begin{pmatrix}
s_i\\
\tau_i\\
p_{i,j}\\
r_j
\end{pmatrix}$. Then we get the following $4$-network of the $t_i$'s decision
\begin{figure}[h!]
\centering
{%
\resizebox*{15cm}{!}{\includegraphics{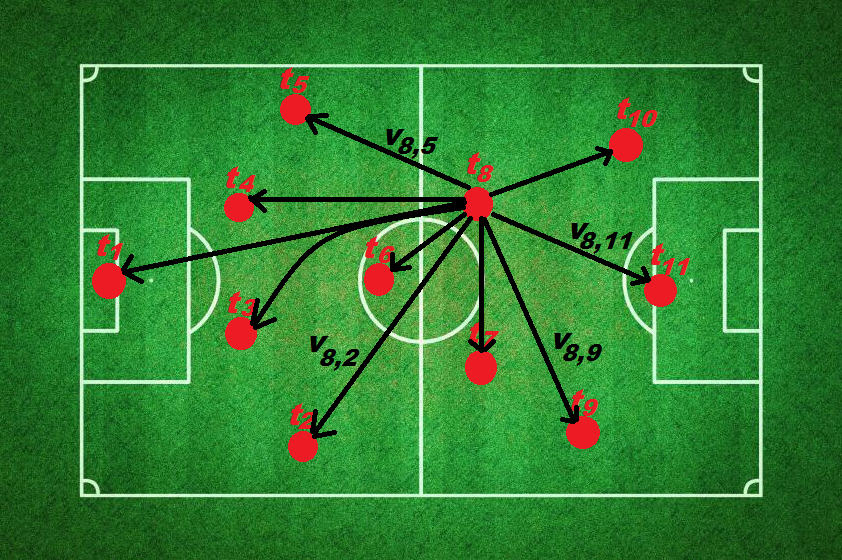}}}\hspace{3pt}
\caption{A $4$-network of the player $t_8$} \label{figure2}
\end{figure}

$$N(i)=\{(\{t_i,t_j\}, v_{i,j} )\mid i\neq j \emph{ and } 1\leq j \leq 11\}.$$
The decision of the player $t_i$ depends on the game style of his team. Then, we determine the decision of the player $t_i$ by the following two functions: 
\begin{enumerate}
\item The game style function $S: (p,r) \mapsto S(p,r)$. By this function, we specify the value of importance for each parameter $p$ and $r$ in the game style of the team. Notice that $0\leq r \leq 10$ and $0 \leq p \leq 1$. Then, to normalize these two parameters, we multiply $p$ by $10$, so we get $0 \leq 10p\leq 10$. A linear game-style function has the following form $$S(p,r)= x \times 10p + y \times r,$$ for some integers $x$ and $y$. Then the importance of the parameter $p$ in this function can be measured by $imp(p)=\frac{x}{x+y}$, and the importance of $r$ can be measured by $imp(r)=\frac{y}{x+y}$. So we have two cases:
\begin{itemize}
\item[$\bullet$] \textbf{Case 1:} if $x \geq y$, then the parameter $p$ is more important than the parameter $r$. This means that in the game style designed by this function, the team aims to possess the ball more than create the risk, we call it a \textit{possession style}.
\item[$\bullet$] \textbf{Case 2:} if $x \leq y$, then the parameter $r$ is more important than the parameter $p$. This means that in the game style designed by this function, the team aims to create the risk more than possessing the ball, we call it a \textit{direct style}.
\end{itemize}    
\item The decision function $D$ that we define by $$ D(t_i) = \left\{\begin{array}{ll}
\emph{shot}&\emph{if } s_i \geq \tilde{s}, \\
\emph{pass to }t_j &\emph{if } s_i < \tilde{s} \emph{ and } S( p_{i,j}, r_j) = \max \limits_{1 \leq k \leq 11 \atop i \neq k} (S( p_{i,k}, r_k)),
\end{array}
\right.$$
for a fixed real number $\tilde{s} \in [0,1]$. If the probability $s_i$ to score is larger than $\tilde{s}$, then the player $t_i$ should shoot the ball. However, if the probability $s_i$ to score is less than $\tilde{s}$, then the player $t_i$ should pass the ball to his teammate $t_j$ who has the optimal game style function. 
\end{enumerate}

Suppose now that the player $t_i$ passes the ball to the player $t_j$ where $i\neq j$. Then, the $4$-network $N(i)$ will be transformed to the $4$-network $N(j)$. Let us consider $N_0$ to be the initial $4$-network, then we get a sequence $(N_i)_{0\leq i\leq n}$ of $4$-networks. So, we should analyze two things about these sequences: 
\begin{itemize}
\item[$\bullet$] \textit{The $s$-efficiency of the sequence}, for $s\in [0,1]$: we aim by this to show if the sequence $(N_i)_{0\leq i\leq n}$ of $4$-networks will give us a certain $4$-network $N_i$ such that the probability of scoring $s_i\geq s$. For example, a $1$-efficient sequence will provide a $4$-network $N_i$ having $s_i=1$, that is this sequence is perfectly efficient since it finishes by scoring a goal. However, this is not really realistic, so if we get for example only a $(0.8)$-efficient sequence, then this will be a very good job since we get an attempt to score with a probability higher than $0.8$ (i.e. $80\%$).
\item[$\bullet$] \textit{The $p$-security of the sequence}, for $p\in [0,1]$: here we aim to show that for any $4$-network $N_i$ of the sequence $(N_i)_{0\leq i\leq n}$, the probability $p_i$ that the pass of the player associated to $N_i$ will be received by the next player, verifies $p_i \geq p$. This allows us to measure the capability of possessing the ball. For example, a $1$-secure sequence means that each pass has a probability $p_i=1$ to be successful, that is the sequence is perfectly secure, and the adversary cannot retrieve the ball. But this is too optimistic, so practically it is enough to construct some $(0.8)$-secure sequences for example, or $(0.6)$-secure sequences if we need more courage.     
\end{itemize}
A successful tactical work aims to construct $s$-efficient and $p$-secure sequences of $4$-networks for some acceptable indices $s$ and $p$. Sometimes, it is not very difficult to construct a $p$-secure sequence but with low efficiency. In this case, if the team need to score more goals, this system will not be relevant. And sometimes, we can construct an $s$-efficient sequence, but with low security. In this case, if we need possession of the ball, we risk losing it very quickly. So, the challenge is to get an optimal pair $(s,p)$ in order to be both secure and efficient.  

\section{Conclusion}
In this paper, we gave a network model for the decision-making of a player who possesses the ball. We showed that this depends on the game style adopted by the team and its coach. This work may contribute to developing player's decisions and also the coach's tactical plans. 

On the other hand, we believe that it is possible to create a similar network model for the decision-making of players without possession. This could enhance the pressing and defensive plans.


\section*{Disclosure statement}
No conflict of interest was reported by the authors.

\section*{Funding}
The author received no financial support for the research, authorship, and/or publication of this article.


\end{document}